\documentclass[onefignum,onetabnum]{siamart171218}

\usepackage{tikz}
\usetikzlibrary{cd}
\usepackage{lipsum}
\usepackage{amsfonts}
\usepackage{graphicx}
\usepackage{epstopdf}
\usepackage{algorithmic}
\ifpdf
  \DeclareGraphicsExtensions{.eps,.pdf,.png,.jpg}
\else
  \DeclareGraphicsExtensions{.eps}
\fi

\usepackage[top=3.375cm,
            bottom=3.375cm,
            inner=4.275cm,
            outer=4.275cm]{geometry}

\newsiamremark{remark}{Remark}
\newsiamremark{hypothesis}{Hypothesis}
\crefname{hypothesis}{Hypothesis}{Hypotheses}
\newsiamthm{claim}{Claim}
\newsiamthm{assumption}{Assumption}
\newsiamremark{example}{Example}

\headers{Fundamental Solutions}{D. S. Winterrose}

\title{A structure theorem for fundamental solutions \\ of analytic multipliers in \texorpdfstring{$\mathbb{R}^n$}{}}

\author{David Scott Winterrose \thanks{Department of Applied Mathematics and Computer Science, Technical University of Denmark, Kgs. Lyngby, Denmark (\email{dawin@dtu.dk}).}}

\usepackage{amsopn}
\usepackage{mathrsfs}

\ifpdf
\hypersetup{
  pdftitle={Fundamental Solutions},
  pdfauthor={David Winterrose}
}
\fi

\newcommand{\Op}{\textnormal{Op}}
\newcommand{\supp}{\textnormal{supp}}

\begin{document}

\maketitle

\begin{abstract}
Using a version of Hironaka's resolution of singularities for real-analytic functions, any elliptic multiplier $\Op(p)$ of order $d>0$, real-analytic near $p^{-1}(0)$, has a fundamental solution $\mu_0$. We give an integral representation of $\mu_0$ in terms of the resolutions supplied by Hironaka's theorem.
This $\mu_0$ is weakly approximated in $H^t_{\textnormal{loc}}(\mathbb{R}^n)$ for $t<d-\frac{n}{2}$
by a sequence from a Paley-Wiener space.
In special cases of global symmetry, the obtained integral representation can be made fully explicit, and we use this to compute fundamental solutions for two non-polynomial symbols.
\end{abstract}

\begin{keywords}
  Fundamental Solutions; Pseudo-Differential Equations; PDE.
\end{keywords}

\begin{AMS}
  35A08, 35E05, 35C05, 35A17
\end{AMS}

\section{Introduction}
Let $ \mathcal{S}'(\mathbb{R}^n)$ denote the space of tempered distributions on $ \mathbb{R}^n$.
A fundamental solution of $\Op(p) = \mathcal{F}^{-1}p\mathcal{F}$ is a $\mu_0 \in \mathcal{S}'(\mathbb{R}^n)$ such that
\begin{align*}
\Op(p)\mu_0
= \delta_0
\quad
\textnormal{in} 
\quad 
\mathcal{S}'(\mathbb{R}^n),
\end{align*}
where $\delta_0$ is the unit measure at $0$, $\mathcal{F}$ is the Fourier transform, and $p$ is the symbol.
The study of these is classical, and most results are recorded in standard texts \cite{HI, HII}.
The H\"ormander-Lojasiewicz theorem \cite{hormander1958solution, lojasiewiez1959solution} ensures existence when $p$ is a polynomial, and provides a way to construct a $\mu_0$, at least in principle, explicitly from the symbol.
But the situation becomes nebulous when $p$ is not a polynomial or globally smooth. We address this problem when $p^{-1}(0)$ is compact and $p$ is real-analytic near $p^{-1}(0)$, and obtain an integral representation that clearly shows the structure of $\mu_0$.  \\

In order to do so, we use a variant of Hironaka's resolution of singularities \cite{atiyah1970}.
Generally, the local charts that are supplied by Hironaka's theorem are unknowable, but they allow us to build an integral representation out of the geometry of $p^{-1}(0)$.
Occasionally, a diffeomorphism that brings $p$ into a resolved form can replace them, and such a diffeomorphism can often be constructed
when there is global symmetry. 
Important examples with this property include any sum of powers of the Laplacian, or any sum of powers of certain elliptic self-adjoint second-order differential operators.
However, the novel and most interesting case for us here is when $p$ is \textit{not} a polynomial, and $p$ is not necessarily globally smooth, but with dimension $n>1$ and order $d>0$. We give examples showing the utility of this approach.\\

A lot of research has been devoted to the construction of explicit representations. See e.g. Ortner and Wagner \cite{ortnerwagner1997} and Camus \cite{camus2006fundamental, camus2008fundamental} for a broad class of operators.
Usually, it is very difficult to find explicit representations of fundamental solutions, and the study is often focused on 
a particular operator of fixed order and dimension. In the case of general homogeneous elliptic and some types of non-elliptic operators,   
Camus \cite{camus2006fundamental, camus2008fundamental} obtained explicit representations valid for any number of dimensions. Apart from the base practical value of constructing general solutions via convolution,
an explicit form may find application in proofs of mapping properties of its operator. See e.g. Rabier~\cite{rabier2015}, where the solution obtained in \cite{camus2006fundamental}, implicit in \cite{HI}, is used.

\newpage
\section{Notation} 
%A symbol, in this article, is a 
Let $p\in C(\mathbb{R}^n)$ be real-analytic in a neighbourhood of $p^{-1}(0) \neq  \emptyset$.
It must be
smooth outside an open ball $B(0,R)$ centered at $0$ of some radius $R>0$.
Putting $\langle \xi \rangle = (1+ |\xi|^2)^\frac{1}{2}$ for $\xi \in \mathbb{R}^n$, it must satisfy
\begin{align*}
\sup_{\xi\in \mathbb{R}^n\setminus B(0,R)} \langle \xi \rangle^{-d} | \partial_\xi^\alpha p(\xi) | < \infty \quad \textnormal{for all} \quad \alpha \in \mathbb{N}_0^n,
\end{align*}
and the ellipticity constraint
\begin{align*}
\inf_{\xi \in \mathbb{R}^n\setminus B(0,R)} \langle \xi \rangle^{-d}|p(\xi)| > 0.
\end{align*}

\begin{definition}[The Paley-Wiener spaces]
Let $K \subset \mathbb{R}^n$ be compact and convex. Define $\textnormal{PW}^d_K(\mathbb{R}^n)$ to be the space of entire functions $u$ satisfying
\begin{align*} 
\sup_{x \in \mathbb{C}^n}\exp\Big(-\sup_{\xi \in K} \textnormal{Im}(x) \cdot \xi \Big) \langle x \rangle^{-d}|u(x)| < \infty.
\end{align*}
If $\{ K_j \}_{j=1}^\infty$ is an exhaustion of $\mathbb{R}^n$ by compact convex sets, we put
\begin{align*}
\textnormal{PW}^d(\mathbb{R}^n) &= \cup_{j=1}^\infty  \textnormal{PW}^d_{K_j}(\mathbb{R}^n).
\end{align*}
Moreover, we put $\textnormal{PW}^{-\infty}(\mathbb{R}^n) = \cap_{d\in \mathbb{Z}} \, \textnormal{PW}^{d}(\mathbb{R}^n)$ and $\textnormal{PW}^{+\infty}(\mathbb{R}^n) = \cup_{d\in \mathbb{Z}} \, \textnormal{PW}^{d}(\mathbb{R}^n)$.
\end{definition}

These spaces are related to $ \mathcal{E}'(\mathbb{R}^n)$, the compactly supported distributions on $\mathbb{R}^n$. 
We write $H^{t}_\textnormal{loc}(\mathbb{R}^n)$ for the Frechet space of distributions locally belonging to $H^{t}(\mathbb{R}^n)$,
and $H^{-t}_\textnormal{comp}(\mathbb{R}^n)$ for its dual space of compactly supported distributions in $H^{-t}(\mathbb{R}^n)$. %, which is equipped with the usual inductive limit topology, in the same way as $C^\infty_0(\mathbb{R}^n)$.
Finally, $ \mathcal{S}(\mathbb{R}^n)$ is the Schwartz space, $u\in \mathcal{S}(\mathbb{R}^n)$ decays faster than any polynomial,
and we put $\mathcal{F}u(\xi) = \int_{\mathbb{R}^n} e^{-ix\cdot \xi} u(x) \, dx$ and $\mathcal{F}^{-1}u(x) = \frac{1}{(2\pi)^n}\mathcal{F}u(-x)$ for $x,\xi\in \mathbb{R}^n$. 

\section{Solution Operator} 
Our main tool is Hironaka's resolution of singularities. 
Often it is stated abstractly \cite{Bjork}, but we need a local embedded version.

\begin{theorem}[Local embedded version of Hironaka's theorem. From Atiyah \cite{atiyah1970}] 
Let $U\subset \mathbb{R}^n $ be an open neighbourhood of $0$, and let $f$ be a function $0 \not\equiv  f \in C^\omega(U)$. 
%Then there exists an open subset $V\subset U$ containing $0$, a real-analytic manifold $M$, and a proper real-analytic map 
Then there is an open $0 \in V\subset U$, a real-analytic manifold $M$, and a map 
\begin{align*}
    \Psi : M \to V.
\end{align*}
It has the following properties:
\medskip
\begin{enumerate}

\item $\Psi : M \to V$ is proper and real-analytic.
\medskip
\item $\Psi : M \setminus (f\circ\Phi)^{-1}(0) \to V \setminus f^{-1}(0) $ is a real-analytic diffeomorphism.
\medskip
\item $(f\circ\Psi)^{-1}(0)$ is a hypersurface in $M$ with normal crossings.

\end{enumerate}
\medskip
\end{theorem}

As $p$ is real-analytic near its zero-set, it is compact with Lebesgue measure zero. The resolution theorem implies that $p$ can be written locally in normal crossing form. 
Fix an open cover $\{V_j\}_{j=1}^N$ of $p^{-1}(0)$ and open $\{U_j\}_{j=1}^N$ such that
\medskip
\begin{enumerate}

\item $\Psi_j : U_j\setminus (p\circ \Psi_j)^{-1}(0) \to V_j\setminus p^{-1}(0)$ is a real-analytic diffeomorphism,
\medskip
\item $(p\circ \Psi_j)(x) = c_j(x) x^{\alpha_j}$ for all $x\in  U_j$ for some $\alpha_j \in \mathbb{N}_0^n$, 

\end{enumerate}
\medskip
where each $c_j$ is a complex-valued, but nowhere zero, real-analytic function on $U_j$. Also, we put $m=\max\{\alpha_j\}_{j=1}^N$.

\newpage
\begin{theorem}\label{thm:mainthmsolop}
Let $\{\chi_j \}_{j=1}^N$ be any partition of unity subordinate to $\{V_j\}_{j=1}^N$.
There is a fundamental solution $\mu_0$, smooth away from $x = 0$, of the form
\begin{align*}
\mu_0(x) = \mathcal{F}^{-1}\Big(\frac{\chi}{p}\Big)(x) +\sum_{j=1}^N\int_{\mathbb{R}^n} I_j(z) \partial_z^{\alpha_j} \Big[ e^{ix\cdot \Psi_j(z)} \frac{(\chi_j \circ \Psi_j)(z)}{c_j(z)} |\det d\Psi_j(z)| \Big] \, dz, 
\end{align*}
where $\chi = 1-\sum_{j=1}^N \chi_j$, and the $I_j$ are given a.e. by
\begin{align*}
I_j(z) =  \frac{1}{(2\pi)^n} \prod_{\alpha_{j,k} \neq 0} \frac{-\ln |z_k| }{(\alpha_{j,k}-1)!}.
\end{align*}
It is weakly approximated in $H^t_{\textnormal{loc}}(\mathbb{R}^n)$ for $t< d- \frac{n}{2}$ by a sequence in $\textnormal{PW}^{m}(\mathbb{R}^n)$.
Finally, if $p^{-1}(0)$ is embedded, any
$v \in \ker \Op(p) \subset \textnormal{PW}^\infty(\mathbb{R}^n)$ is of the form
\begin{align*}
    v(x) 
    =\sum_{j=1}^N  \sum_{k\leq k_j-1}  \Big\langle  (\Psi_j^*)^{-1} (u_{j,k} \otimes \partial_{z_n}^k \delta_0)(\xi),  e^{ix\cdot \xi}  \Big\rangle,
\end{align*}
where $u_{j,k} \in \mathcal{E}'(U_j^0)$ are supported in the $z_n = 0$ slice $U_j^0 = \{ z \in \mathbb{R}^{n-1} \, | \, (z,0) \in U_j \}$, and each $\Psi_j$ is arranged so that 
\begin{align*}
(p \circ \Psi_j) (x) = c_j(x) x^{k_j}_n.
%\quad
%\textnormal{for all}
%\quad
%x\in  U_j.
\end{align*}
In this case, all other fundamental solutions differ from $\mu_0$ by such a $v$ .
\end{theorem}

The first step is to prove a lemma about principal value integrals with log kernel. It is used here in a way similar to Bj\"ork \cite[Chapter 6, Theorem~1.5]{Bjork}.

\begin{lemma} \label{lmm:logres}
Let $\psi \in C^\infty(\mathbb{R})$ be either rapidly decaying or compactly supported. Then, for any $k \in \mathbb{N}$, we have
\begin{align*} 
\int_{-\infty}^\infty \psi(r) \, dr
=
\frac{-1}{(k-1)!}
\int_{-\infty}^\infty  \ln(|r|) \frac{d^k}{dr^k} \Big( r^{k} \psi(r) \Big) \, dr.
\end{align*}
\end{lemma}
\begin{proof}
The proof of this is a routine exercise in repeated integration by parts. Observe that we can write $\psi(r) = r^{-k}r^k \psi(r)$, and 
\begin{align*} 
\int_{-\infty}^\infty \psi(r) \, dr
%&=
%\int_{-\infty}^\infty  \frac{d}{dr} \Big( \frac{-1}{k-1} r^{-k+1}  \Big) r^{k} \psi(r) \, dr \\
&=
\frac{-1}{k-1} r^1 \psi(r) \Big|_{-\infty}^\infty + \frac{1}{k-1}  \int_{-\infty}^\infty  r^{-k+1}  \frac{d}{dr} \Big( r^{k} \psi(r) \Big) \, dr \\
&\cdots \\
&=
\frac{-1}{(k-1)!} r^{k-1}\psi(r) \Big|_{-\infty}^\infty + \frac{1}{(k-1)!}  \int_{-\infty}^\infty  r^{-1}  \frac{d^{k-1}}{dr^{k-1}} \Big( r^{k} \psi(r) \Big) \, dr \\
&=
\frac{-1}{(k-1)!}
\int_{-\infty}^\infty  \ln(|r|) \frac{d^k}{dr^k} \Big( r^{k} \psi(r) \Big) \, dr,
\end{align*}
where all boundary terms at $0$ in the final integration vanish, because
\begin{align*}
   \lim_{r \to 0\pm} \ln(|r|) \frac{d^{k-1}}{dr^{k-1}} \Big( r^{k} \psi(r) \Big) = 0,
\end{align*}
and boundary terms at $\pm \infty$ vanish by the hypothesis on $\psi$. 
\end{proof}

\newpage
\begin{lemma} \label{lmm:solop}
Define $Q : \mathcal{S}(\mathbb{R}^n) \to C^\infty(\mathbb{R}^n)$  by
\begin{align*} %\label{eq:Qj}
Q v(x)
= \sum_{j=1}^N\int_{\mathbb{R}^n} I_j(z) \partial_z^{\alpha_j} \Big[ e^{ix\cdot \Psi_j(z)} \frac{(\chi_j\mathcal{F}v) \circ \Psi_j (z)}{c_j(z)} |\det d\Psi_j(z)| \Big] \, dz.
\end{align*}
Then $P = \Op(\frac{\chi}{p}) + Q $ satisfies both 
\begin{align*}
    \Op(p)Pv = v
    \quad
    \textnormal{and}
    \quad
    P\Op(p)v = v
    \quad
    \textnormal{for all}
    \quad
    v \in \mathcal{S}(\mathbb{R}^n).
\end{align*}
\end{lemma}
\begin{proof}
The proof is an application of \cref{lmm:logres} and the Fubini-Tonelli theorem.
Let $\psi \in \mathcal{S}(\mathbb{R}^n)$. Using \cref{lmm:logres} coordinate-wise, we compute
\begin{align*}
\langle \Op(p)Q v, \psi \rangle 
&=
\int_{\mathbb{R}^n} I_j(z) \partial_z^{\alpha_j} \Big[ \langle \Op(p) e^{i(\cdot)\cdot \Psi_j(z)} , \psi \rangle \frac{(\chi_j\mathcal{F}v) \circ \Psi_j (z)}{c_j(z)} |\det d\Psi_j(z)| \Big] \, dz \\
&=
\int_{\mathbb{R}^n} I_j(z) \partial_z^{\alpha_j} \Big[ z^{\alpha_j} \langle e^{i(\cdot)\cdot \Psi_j(z)} , \psi \rangle (\chi_j\mathcal{F}v) \circ \Psi_j (z) |\det d\Psi_j(z)| \Big] \, dz  \\
&=
 \frac{1}{(2\pi)^n}\int_{\mathbb{R}^n \setminus (p\circ \Psi_j)^{-1}(0)}  \langle e^{i(\cdot)\cdot \Psi_j(z)} , \psi \rangle (\chi_j\mathcal{F}v) \circ \Psi_j (z) |\det d\Psi_j(z)|  \, dz \\
&=  
\frac{1}{(2\pi)^n}\int_{\mathbb{R}^n \setminus p^{-1}(0)}  \langle e^{i(\cdot)\cdot \xi } , \psi \rangle (\chi_j\mathcal{F}v) (\xi) \, d\xi \\
&= \langle\Op(\chi_j)v, \psi \rangle,
\end{align*}
and we then get 
\begin{align*}
\langle \Op(p)Q v, \psi \rangle
=
\sum_{j=1}^N \langle\Op(\chi_j)v, \psi \rangle = \langle v - \Op(\chi)v, \psi \rangle.
\end{align*}
Point-wise in $x\in \mathbb{R}^n$, we compute
\begin{align*}
Q\Op(p)v(x)
&=
\int_{\mathbb{R}^n} I_j(z) \partial_z^{\alpha_j} \Big[ e^{i x \cdot \Psi_j(z)}  \frac{(\chi_j\mathcal{F}\Op(p)v) \circ \Psi_j (z)}{c_j(z)} |\det d\Psi_j(z)| \Big] \, dz \\
&=
\int_{\mathbb{R}^n} I_j(z) \partial_z^{\alpha_j} \Big[ z^{\alpha_j}  e^{ix\cdot \Psi_j(z)} (\chi_j\mathcal{F}v) \circ \Psi_j (z) |\det d\Psi_j(z)| \Big] \, dz  \\
&=
 \frac{1}{(2\pi)^n}\int_{\mathbb{R}^n \setminus (p\circ \Psi_j)^{-1}(0)}  e^{i x\cdot \Psi_j(z)} (\chi_j\mathcal{F}v) \circ \Psi_j (z) |\det d\Psi_j(z)|  \, dz \\
&=
\frac{1}{(2\pi)^n}\int_{\mathbb{R}^n \setminus p^{-1}(0)}  e^{ix\cdot \xi } (\chi_j\mathcal{F}v) (\xi) \, d\xi \\
&= \Op(\chi_j)v(x),
\end{align*}
which shows that
\begin{align*}
Q\Op(p)v = \sum_{j=1}^N \Op(\chi_j)v = v - \Op(\chi)v.
\end{align*}
Note that the properties of $\Psi_j$ ensure that all the above integrals are well-defined. The determinant of $d\Psi_j$ on each component of $U_j \setminus (p\circ \Psi_j)^{-1}(0)$ never becomes zero.
This completes the proof.
\end{proof}

%In particular, the above lemma implies that
%\begin{align*} \label{eq:paradiff}
%\Op\Big(\frac{\chi}{p}\Big)v = P\Op(p)\Op\Big(\frac{\chi}{p}\Big)v  = Pv - P\Op(1-\chi)v.
%\end{align*}
%and in particular,
%\begin{align*}
%Qv = P\Op(1-\chi)v.
%\end{align*}

\newpage
\begin{lemma} \label{lmm:contQj}
$P : \mathcal{S}(\mathbb{R}^n) \to C^\infty(\mathbb{R}^n)$ is continuous.
\end{lemma}
\begin{proof}
The proof is just estimating $C^\infty(\mathbb{R}^n)$ semi-norms of $Q$ in those of $\mathcal{S}(\mathbb{R}^n)$.
By the chain rule, if $v\in \mathcal{S}(\mathbb{R}^n)$, we have
\begin{align*} 
\partial_{z_k}(\mathcal{F}v \circ \Psi_j) 
=
-i(\mathcal{F}(x_1 v)\circ \Psi_j, \cdots,\mathcal{F}(x_n v)\circ \Psi_j) \cdot \partial_{z_k}\Psi_j.
\end{align*}
Using the Leibniz rule, we get for any $\alpha\in \mathbb{N}^n_0$ some $C_\alpha',C_\alpha > 0$ such that
\begin{align*} 
|\partial_{x}^\alpha Q v(x)|
&\leq
 C_\alpha'
 \sum_{j=1}^N\int_{\supp (\chi_j)} |I_j(z)| \sum_{\beta \leq \alpha_j} \langle x\rangle^{\alpha_j - \beta}
 |\partial_z^{\beta}  (\mathcal{F}v \circ \Psi_j) (z)| \, dz\\
&\leq C_\alpha \sum_{j=1}^N
\Big( \int_{\supp (\chi_j)} |I_j(z)| \, dz \Big) \langle x \rangle^{m} \max_{|\beta|\leq m}\sup_{z\in \mathbb{R}^n} | \mathcal{F}(x^\beta v)(z) |,
\end{align*}
which by the continuity of $\mathcal{F} : \mathcal{S}(\mathbb{R}^n) \to \mathcal{S}(\mathbb{R}^n)$ implies the lemma.
\end{proof}

\begin{lemma} \label{lmm:solmapping}
$P : \textnormal{PW}^{-\infty}(\mathbb{R}^n) \to \textnormal{PW}^m(\mathbb{R}^n)$ is well-defined.
\end{lemma}
\begin{proof}
Because each map $\Psi_j$ is proper, each $( \chi_j \mathcal{F}v )\circ\Psi_j$ is compactly supported. 
By the well-known~\cite[Paley-Wiener-Schwartz Theorem~7.3.1]{HI}, $\Op(\frac{\chi}{p})v \in \textnormal{PW}^{-\infty}(\mathbb{R}^n)$.
A simple estimate gives $C',C>0$ such that
\begin{align*}
    |Qv(x)| 
    &\leq C'
    \sum_{j=1}^N\int_{\supp (\chi_j)} |I_j(z)| \sum_{\beta \leq \alpha_j} |\partial_z^{\beta} [ e^{ix\cdot \Psi_j(z)} ] | \, dz
    \\
    &\leq 
    C\sum_{j=1}^N
\Big( \int_{\supp (\chi_j)} |I_j(z)| \, dz \Big) \langle x \rangle^{m}
\exp\Big(\sup_{\xi \in K} \textnormal{Im}(x)\cdot \xi\Big),
%\\  &\leq  C \langle x \rangle^{m} \exp\Big( \sup_{\xi \in K} \textnormal{Im}(x)\cdot \xi \Big),
\end{align*}
%where $K$ is suitable a compact and convex set so large that
%$-\cup_{j=1}^N \supp(\chi_j) \subset K$,
where $K$ is a compact and convex set so large that
\begin{align*}
    -\cup_{j=1}^N \supp(\chi_j) \subset K,
\end{align*}
and so $Qv$ is entire with $Qv \in \textnormal{PW}^{m}(\mathbb{R}^n)$.
\end{proof}

Define the reflection map $A$ by $A\psi(x) = \psi(-x)$ for all $x\in \mathbb{R}^n$ on $\psi \in C^\infty(\mathbb{R}^n)$.
It takes $\mathcal{S}(\mathbb{R}^n)$ and $C^\infty_0(\mathbb{R}^n)$ continuously to themselves. 
The transpose of $Q$ is $AQA$.
%and similarly for $\Op(\frac{\chi}{p})$, since $\mathcal{F}A = (2\pi)^{n} A \mathcal{F}^{-1}$ by the definition of $\mathcal{F}$ and $\mathcal{F}^{-1}$.
Using this fact and \cref{lmm:contQj}, we extend $Q$, hence $P$, by duality:
\begin{definition} \label{def:solopdual}
Define  $Q: u\mapsto Qu : \mathcal{E}'(\mathbb{R}^n) \to \mathcal{S}'(\mathbb{R}^n)$ by
\begin{align*}
\langle Qu, \psi \rangle = \langle u, A Q A\psi \rangle
\quad
\textnormal{for all}
\quad
\psi \in \mathcal{S}(\mathbb{R}^n).
\end{align*}
\end{definition}

\begin{lemma} \label{lmm:solopdistrib}
$\Op(p)Ps = s$ holds for any $s\in \mathcal{E}'(\mathbb{R}^n)$.
\end{lemma}
\begin{proof} 
%Since $\Op(p) Ps =\Op(\chi)s $, we must show that $\Op(p)P s = \Op(1-\chi)s$.
By \cref{lmm:solop}, if $\psi \in \mathcal{S}(\mathbb{R}^n)$, we have
\begin{align*}
\langle \Op(p)P s , \psi \rangle
&=
\langle s, A P A \Op(p)^* \psi \rangle \\
&=
\langle s, A P A \mathcal{F} (p \mathcal{F}^{-1}\psi) \rangle \\
&=
\langle s, A P \Op(p) A\psi \rangle \\
&=
\langle s, A^2\psi \rangle \\
&=
\langle s, \psi \rangle.
\end{align*}
\end{proof}

\newpage
Applying \cref{lmm:solopdistrib}, we get the fundamental solution $\mu_0 = P\delta_0$ for the operator. Using e.g.~\cite[Theorems~5.2 and 7.1]{Shubin2001}, or similar in \cite{GG}, it is smooth in $x \neq 0$, and 
\begin{align*}
     \mu_0 \in H^{t}_{\textnormal{loc}}(\mathbb{R}^n)
    \quad
    \textnormal{if}
    \quad
    t < d-\frac{n}{2}.
\end{align*}

\begin{lemma} 
$\mu_0$ is weakly approximated in $H^t_{\textnormal{loc}}(\mathbb{R}^n)$ by a $\textnormal{PW}^{m}(\mathbb{R}^n)$ sequence.
\end{lemma}
\begin{proof}
Take a bump function $\eta \in C^\infty_0(\mathbb{R}^n)$ such that $\eta(x) = 1$ holds for all $|x|<1$. Put $\eta_k(x) = \eta(\frac{x}{k})$ for all $x \in \mathbb{R}^n$ and $k\in \mathbb{N}$.
By \cref{lmm:solmapping}, $P\mathcal{F}^{-1}\eta_k \in \textnormal{PW}^{m}(\mathbb{R}^n)$.
Given any $u\in H^{-t}_{\textnormal{comp}}(\mathbb{R}^n)$, then for $k$ large enough, we get
\begin{align*}
    |\langle \mu_0 - P\mathcal{F}^{-1}\eta_k, u \rangle|^2  
    &= 
    \Big|\Big \langle \mathcal{F}^{-1}\Big(\frac{\chi}{p}(1- \eta_k)\Big) , u \Big \rangle\Big|^2 \\
    &= 
    \Big|\Big \langle \frac{\chi}{p}(1- \eta_k) , \mathcal{F}^{-1}u  \Big \rangle\Big|^2 \\
    &\leq
    \Big( \int_{\mathbb{R}^n} \langle \xi \rangle^{2t} \Big|\frac{\chi}{p}(1- \eta_k)\Big|^2 \, d\xi \Big) \Big( \int_{\mathbb{R}^n} \langle \xi \rangle^{-2t} |\mathcal{F}^{-1}u(\xi)|^2 \, d\xi \Big)
\end{align*}
and so $P\mathcal{F}^{-1}\eta_k \to \mu_0$ weakly in $H^t_{\textnormal{loc}}(\mathbb{R}^n)$ as $k\to \infty$.
\end{proof}

\begin{lemma} 
Suppose that $p^{-1}(0)$ is embedded as a real-analytic submanifold. Then $\ker \Op(p)$ consists of functions $v \in  \textnormal{PW}^{\infty}(\mathbb{R}^n)$ of the form
\begin{align*}
    v(x) 
    = \sum_{j=1}^N  \sum_{k\leq k_j-1}  \Big\langle (\Psi_j^*)^{-1} (u_{j,k} \otimes \partial_{z_n}^k \delta_0)(\xi),  e^{ix\cdot \xi}  \Big\rangle,
\end{align*}
where $u_{j,k} \in \mathcal{E}'(U_j^0)$ and $\Psi_j$ are precisely as stated in \cref{thm:mainthmsolop}.
\end{lemma}
\begin{proof}
Observe that $p\mathcal{F}v = 0$ implies $\textnormal{supp}\, \mathcal{F}v \subset p^{-1}(0)$ so that $\mathcal{F}v$ is compact. Again by ~\cite[Theorem~7.3.1]{HI}, $v \in  \textnormal{PW}^{\infty}(\mathbb{R}^n)$. Observe then that
\begin{align*}
0 = \Psi_j^* ( \chi_j p\mathcal{F}v )
= c_j z^{k_j}_n \Psi_j^* (\chi_j \mathcal{F}v),
\end{align*}
and since $c_j$ is never zero, by~\cite[Theorem~2.3.5]{HI}, we must have
\begin{align*}
\Psi_j^* ( \chi_j \mathcal{F}v )
    =
    \sum_{k\leq k_j-1} (2\pi)^n u_{j,k} \otimes \partial_{z_n}^k \delta_0,
\end{align*}
where $u_{j,k} \in \mathcal{E}'(\mathbb{R}^{n-1})$ are some distributions supported inside the $z_n = 0$ slice of $U_j$. It follows that
\begin{align*}
    v(x) 
    &= \mathcal{F}^{-1}\Big(\sum_{j=1}^N \chi_j \mathcal{F}v\Big)(x) \\
    &= \sum_{j=1}^N  \sum_{k\leq k_j-1} \mathcal{F}^{-1} (\Psi_j^{-1})^* \Big((2\pi)^n u_{j,k} \otimes \partial_{z_n}^k \delta_0\Big)(x) \\
    &= \sum_{j=1}^N  \sum_{k\leq k_j-1}  \Big\langle (\Psi_j^*)^{-1} (u_{j,k} \otimes \partial_{z_n}^k \delta_0)(\xi),  e^{ix\cdot \xi}  \Big\rangle.
\end{align*}
\end{proof}

\newpage

\begin{figure}[H]
    \centering
    \begin{tikzpicture}

      \draw[thick, domain=0:2*pi,samples=500] plot ({deg(\x)}:{1+0.25*cos(6*\x r)-0.25*sin(4*\x r) });
      
      %\draw[dashed, domain=0:2*pi,samples=500] (0,0) circle (1);
      \draw[thick, domain=0:2*pi,samples=500] (5,0) circle (1);
      
      \filldraw [white, even odd rule, fill=gray!75, fill opacity=.3,domain=0:2*pi,samples=200,smooth ] plot (xy polar
      cs:angle=\x r,radius= { 1.25+0.25*cos(6*\x r)-0.25*sin(4*\x r) })
      plot (xy polar cs:angle=\x r,radius= { 0.75+0.25*cos(6*\x r)-0.25*sin(4*\x r) });
      
      \filldraw[white, even odd rule, fill=gray!75, fill opacity=.3, samples=200,smooth](5,0) circle (1.25) (5,0) circle (0.75);
      
      %\fill[white, fill=gray!75, fill opacity=.3,samples=200,smooth] (8.75,1.25) rectangle (9.25,-1.25);
      %\draw[thick] (9,1.25)--(9,-1.25);
      
      \draw[thick,->,>=latex] (3,0)--(2,0);
      %\draw[thick,->,>=latex] (7,0)--(8,0);
      
      \node at (0,-1.75)  {$p^{-1}(0)$};
      \node at (5,-1.75)  {$(p\circ\Psi)^{-1}(0)$};
      \node at (2.5,-0.5)  {$\Psi$};
      
      %\draw[->] (7,-1) .. controls (7.5,-1) .. (8,-1.5);
      
    \end{tikzpicture}
    \caption{Deformation of a star-convex zero-set onto a circle.}
\end{figure}

The main \cref{thm:mainthmsolop} is finally obtained by combining the above partial results.
%Unfortunately, in these results $\Psi_j$ can not be directly constructed or obtained easily. In general, it is impossible to obtain them explicitly for any given multiplier symbol. 
Unfortunately, it is impossible to obtain $\Psi_j$ explicitly for any given multiplier symbol. 
But if, for example, $p^{-1}(0)$ is the real-analytic boundary of some star-convex domain, we can replace the charts by a single deformation $\Psi$  of the boundary onto a sphere.
Given $p$, we look for $\Psi$ so that $\Psi^* p$ factorizes.
Our main theorem gives
\begin{align*}
    \mu_0(x) = 
    \mathcal{F}^{-1}\Big(\frac{\chi}{p}\Big)(x) +
    Q\delta_0(x),
\end{align*}
where $\chi$ appropriately suppresses a region surrounding $p^{-1}(0)$ on which $\Psi$ is defined.
%In the following, we give examples to demonstrate this.

\subsection{Sums of powers of \texorpdfstring{$\Delta_g$}{the Laplacian}}
Let $g$ be a positive-definite symmetric matrix.
Consider for $d\in \mathbb{N}$ and $\{c_j\}_{j=0}^{d}\subset \mathbb{C}$ with $c_d = 1$ and $c_0 \neq 0$ the multiplier
\begin{align*}
    \Op(p) = \sum^d_{j= 0} c_j \Delta_g^\frac{j}{2}.
\end{align*}
The symbol $p$ is taken into a polynomial form by the map
\begin{align*}
   \Psi : (0,\infty) \times \mathbb{S}^{n-1} \to \mathbb{R}^{n}\setminus 0 : (r,\omega) \mapsto r g^{-\frac{1}{2}} \omega.
\end{align*}
Pulling back, we find that
\begin{align*}
    (p \circ \Psi)(r,\omega) 
    %= \sum^d_{j=0} c_j r^{j} 
    = c(r) \prod^m_{j=1} (r-r_j)^{m_j},
\end{align*}
where $c$ is a polynomial with no root in $[0,\infty)$, and $r_j>0$ are the positive real roots.
Let $C_{j,k}$ be the unique coefficients in $\prod^m_{j=1} (r-r_j)^{-m_j} = \sum_{j=1}^m \sum^{m_j}_{k=1} C_{j,k} (r-r_j)^{-k}$.
Pick $\chi\in C^\infty(\mathbb{R}^n)$ so that $1-\chi \circ \Psi \in C^\infty_0((0,\infty)\times \mathbb{S}^{n-1})$
is independent of $\omega\in \mathbb{S}^{n-1}$, and
equal to $1$ in a neighbourhood of $\cup_{j=1}^m \{r_j\} \times \mathbb{S}^{n-1}$, all of the spheres of radius $r_j$.
If the multiplicities satisfy $m_j< n$, we have
\begin{align*}
Q\delta_0(x)
= \sum_{j=1}^{m} \sum_{k=1}^{m_j} B_{j,k} \int^\infty_0 \ln |r-r_j| \, \partial_r^k \Big[ (1-\chi\circ\Psi)(r) \frac{r^{\frac{n}{2}}}{c(r)} 
\frac{J_{\frac{n}{2}-1}(r|g^{-\frac{1}{2}}x|)}{|g^{-\frac{1}{2}}x|^{\frac{n}{2}-1}} 
\Big] \, dr, 
\end{align*}
where $J_{\frac{n}{2}-1}$ is the cylindrical Bessel function of order $\frac{n}{2}-1$, and
\begin{align*}
B_{j,k}
=
-\frac{\det \, g^{-\frac{1}{2}}}{(2\pi)^{\frac{n}{2}}(k-1)!} C_{j,k}.
\end{align*}
%\begin{align*}
%\int_{\mathbb{S}^{n-1}} e^{i x\cdot r g^{-\frac{1}{2}}\omega} \, \textnormal{vol}_{\mathbb{S}^{n-1}}(\omega)
%=
%(2\pi)^{\frac{n}{2}} \frac{J_{\frac{n}{2}-1}(|r g^{-\frac{1}{2}} x|)}{|r g^{-\frac{1}{2}} x|^{\frac{n}{2}-1}}.
%\end{align*}

\newpage
\subsection{A perturbation of \texorpdfstring{$\Delta_g$}{the Laplacian}}
Let $\arg\xi$ be the multi-valued argument of $\xi\in \mathbb{R}^2$. 
Consider instead the multiplier symbol 
\begin{align*}
    p : \mathbb{R}^2 \to \mathbb{R} : \xi \mapsto |\xi|_g^2 - \Big(1+a\cos(n \arg\xi) \Big),
\end{align*}
where $n \in \mathbb{N}$ and $a<\frac{1}{2}$. 
It is certainly real-analytic near its star-shaped zero-set. 
This $p$ is taken into normal crossing form by the map
\begin{align*}
    \Psi : (-{\textstyle\frac{1}{2}},{\textstyle\frac{1}{2}}) \times (0, 2\pi) \to \mathbb{R}^2 : (r,\theta) \mapsto \Big(r + 1 + a\cos(n \theta) \Big)^\frac{1}{2} g^{-\frac{1}{2}} (\cos\theta, \sin\theta).
\end{align*}
\begin{figure}[H]
    \centering
    \begin{tikzpicture}

      \filldraw [white, even odd rule, fill=gray!75, fill opacity=.3,domain=0:2*pi,samples=200,smooth ] plot (xy polar
      cs:angle=\x r,radius= { sqrt(1.5+0.25*cos(12*\x r)) })
      plot (xy polar cs:angle=\x r,radius= { sqrt(0.5+0.25*cos(12*\x r)) });
      
      \fill[fill=gray!75, fill opacity=.3,samples=200,smooth] (3.5,1.25) rectangle (4.5,-1.25);
      \draw[thick] (4,1.25)--(4,-1.25);
      
      \draw[thick, domain=0:2*pi,samples=500] plot ({deg(\x)}:{sqrt(1+0.25*cos(12*\x r))});
      \node at (0,-1.75)  {$p^{-1}(0)$};
      \draw[thick,<-,>=latex] (2,0)--(3,0);
      
      \node at (4,-1.75)  {$(p\circ \Psi)^{-1}(0)$};
      \node at (2.5,-0.5)  {$\Psi$};
      \draw[thick,white] (-1.375+2,0)--(-0.625+2,0);
      
    \end{tikzpicture}
    \caption{Covering the zero-set of $p$ except for a point. Here $a=\frac{1}{4}$ and $n=12$.}
\end{figure}

\noindent  It is clear that $\Psi$ is a diffeomorphism onto its image, not covering the whole zero-set.
But a representation using only $\Psi$ is still possible, because it misses just a single point.
Pulling back, we find that
\begin{align*}
    ( p \circ \Psi )(r,\theta) = r,
\end{align*}
and we compute
\begin{align*}
    \det d\Psi(r, \theta) 
    %= \frac{1}{|\Psi(r,\theta)|^{2}}
    %\det
    %\begin{bmatrix}
    %\frac{1}{2} \cos\theta & -\frac{1}{2}  a n\sin(n \theta) \cos\theta - |\Psi(r,\theta)|^{2} \sin \theta \\
    %\frac{1}{2}  \sin\theta
    %& -\frac{1}{2}  a n\sin(n \theta) \sin\theta + |\Psi(r,\theta)|^{2} \cos \theta
    %\end{bmatrix}
    =
    \frac{1}{2} \det g^{-\frac{1}{2}}.
\end{align*}
Pick $\chi$ such that $1-\chi \circ \Psi \in C^\infty((-{\textstyle\frac{1}{2}},{\textstyle\frac{1}{2}})\times(0,2\pi))$
does not depend on $\theta\in(0,2\pi)$, and is compactly supported in $(-{\textstyle\frac{1}{4}},{\textstyle\frac{1}{4}})$ and
equal to $1$ in a neighbourhood of $r=0$.
We tacitly extend $\chi$ by one to all of $\mathbb{R}^2$.
In that case, we have
\begin{align*}
    Q\delta_0(x) = - \frac{\det g^{-\frac{1}{2}}}{8\pi^2} \int_{-\frac{1}{2}}^{\frac{1}{2}}\ln |r|\, \partial_r \Big[  (1-\chi \circ \Psi)(r) \int_{0}^{2\pi}e^{ix\cdot \Psi(r,\theta)} \, d\theta  \Big] \, dr,
\end{align*}
and $\ker \Op(p)$ consists of $v$ of the form
\begin{align*}
    v(x) 
    = \Big\langle u(\theta),  e^{ix\cdot \Psi(0,\theta)} \Big\rangle,
\end{align*}
where $u \in \mathcal{D}'(\mathbb{R}/2\pi \mathbb{Z})$ is a distribution on the space of $2\pi$-periodic smooth functions. We could replace $|\xi|^2_g$ in $p$ by any integer power of $|\xi|^2_g$ and still get a similar result, provided that we adjust the fractional power $\frac{1}{2}$ in $\Psi$ in accordance with this change.
A similar technique can be applied to sums of powers of such multipliers too. %arbitrary sums of powers of the Laplacian \cite{KaramehmedovicWinterrose2019}.

\newpage
\section{Acknowledgements}
The author wishes to thank the anonymous reviewer. Their useful comments and suggestions helped to improve this paper.

\bibliographystyle{siamplain}
\bibliography{references}
\end{document}